# Trajectory Convergence from Coordinate-wise Decrease of General Energy Functions

Julien M. Hendrickx [a], Balázs Gerencsér [b]

[a] *ICTEAM Institute, UCLouvain, Louvain-la-Neuve, Belgium*

[b] *Alfréd Rényi Institute of Mathematics, Budapest, Hungary,*
*Eötvös Loránd University, Department of Probability and Statistics, Budapest, Hungary*

**Abstract**

We consider arbitrary trajectories subject to a coordinate-wise energy decrease: the sign of the derivative of each entry is never the same as that of the corresponding entry of the gradient of some energy function. We show that this simple condition guarantees convergence to a point, to the minimum of the energy functions, or to a set where its Hessian has very specific properties. This extends and strengthens recent results that were restricted to convex quadratic energy functions. We demonstrate the application of our result by using it to prove the convergence of a class of multi-agent systems subject to multiple uncertainties.

## 1 Introduction

We consider the convergence properties of a trajectory $y : \mathbb{R}^+ \to \mathbb{R}^n$, $t \mapsto y(t)$ whose evolution is constrained by an energy function $V : \mathbb{R}^n \to \mathbb{R}$, $x \mapsto V(x)$ via the set of inequalities

$$\dot{y}_i \frac{\partial V}{\partial x_i}|_{y(t)} \leq 0, \qquad \forall t \geq 0, i = 1, \ldots, n, \qquad (1)$$

i.e. the derivative of a coordinate $y_i$ of $y$ and the corresponding coordinate of the gradient $\nabla V$ of $V$ at $y$ always have opposite sign if they are both nonzero. We note that we will always use the letter $y$ when referring to a trajectory or its accumulation points, and the letter $x$ for the points in the ambient space $\mathbb{R}^n$.

For energy functions of the form $V = x^T Q x$, with $Q$ positive semi-definite, condition (1) was shown in our recent work [9] to be often sufficient for convergence of $y$. This was motivated by a platoon cooperative control application involving dead-zone control and bounded arbitrary disturbance. It allowed in particular solving a conjecture on the convergence of such systems [5], and a related problem of consensus under bounded disturbance [2,6]. In this work, we extend these results to general functions $V$, characterize more precisely the alternative long-term behavior of the trajectory when convergence is not guaranteed, and explore the tightness of our conditions. We demonstrate the application of our result on a generalization of the platoon application from [9].

We stress that the trajectory $y$ is not assumed to be generated by a vector field or a system of differential equations. It can be completely arbitrary provided it satisfies the constraints (1). By contrast, a large proportion of the convergence results based on decrease of energy functions rely on variations or extensions of Lyapunov-Kraskowski-LaSalle Theorems [10], and typically assume that trajectories follow some ordinary differential equation such as $\dot{y}(t) = f(y(t), t)$ or $\dot{y}(t) = f(y(t))$ for an $f$ satisfying some (uniform) continuity conditions [1,4]. For example, LaSalle theorem guarantees (under some conditions) the convergence of $\dot{y} = f(y)$ to an invariant set, but not to a single point, provided $f(x)^T \nabla V(x) \leq 0$ everywhere [11]. Convergence to 0 can then be guaranteed under the additional assumption that $\frac{d}{dt}V(y(t)) = f(y(t))^T \nabla V(y(t))$ is not uniformly zero *along any trajectory* other than that staying at 0 [14]. For more detail on various cases of unforced systems, we refer the reader

---

to [12] as a starting point.

Vector fields $f(x,t)$ over the state space may not be naturally available in systems whose evolution is driven by external elements. Think of discrete communications in cyber-physical applications, systems designed to be robust to adversarial input signals that could depend on the trajectory and its history, or systems involving some random decisions (though more complex descriptions may be available, see [7,15]). Similarly, many modern control laws are not easily described by a continuous field $f$, think e.g. of event-triggered or self-triggered mechanisms [8,13]. Hence it is desirable to have results guaranteeing the convergence of a *single trajectory based on properties satisfied along that specific trajectory* without assuming or constructing a corresponding vector field, nor speculating about the properties of potential other trajectories. Currently available results for single trajectories require a sufficiently negative decrease, e.g., $\frac{d}{dt}V(y(t)) \leq -\lambda V(y(t))$ for some positive $\lambda$, which allows guaranteeing convergence to the minimum of $V$ at a certain rate, see again [10]. This precludes their use when no such uniform condition can be guaranteed, or for situations where the rate of convergence cannot be known, which could happen for example if parts of the system can occasionally pause. On the other hand, simply requiring $\frac{d}{dt}V(y(t)) < 0$ does not imply convergence, as can be verified on the simple two-dimensional example $y_1(t) = (1+e^{-t})\cos t$ and $y_2(t) = (1+e^{-t})\sin t$ with $V(x) = ||x||_2$. As a source of intuition at a very informal level, one could say that the condition $\frac{d}{dt}V(y(t)) < 0$ implies the decrease and convergence of the energy $V$ along the trajectory, but allows for persistent significant energy transfer between the different coordinates. By contrast, our condition (1) forces the decrease of energy on every coordinate. This remains at the level of intuition though, as the energy $V$ in general cannot be separated along the different coordinates.

Our paper is organized as follows: We state our main convergence result in Section 2, together with convenient corollaries specializing it, including for convex functions $V$. We study its tightness in Section 3 with an example showing that the alternative to convergence cannot simply be removed, and demonstrate in Section 4 its application to analyze a multi-agent problem. The main proof is presented in Section 5, together with the intuition on how the elements are built together. We draw conclusions and discuss potential continuations and open problems in Section 6.

## 2 Main results

We first present our most general result with minimal assumptions, thus allowing for most possibilities for the asymptotic behavior.

**Theorem 1** *Let $V : \mathbb{R}^n \to \mathbb{R}$ be a twice differentiable function with a locally Lipschitz Hessian and $y : \mathbb{R}^+ \to \mathbb{R}^n$ a trajectory that is absolutely continuous, also implying that $\dot{y}(t)$ exists almost everywhere. Suppose that where it exists,*

$$\dot{y}_i(t) \frac{\partial V}{\partial x_i}|_{y(t)} \leq 0 \qquad \forall t \geq 0, i = 1, \ldots, n. \quad (2)$$

*Then, at least one of the following conditions holds:*

*(a) $y$ converges;*
*(b) for every accumulation point $\bar{y}$, there exists a nontrivial vector $v$ in the kernel of the Hessian $\nabla^2 V(\bar{y})$, whose entries $v_i$ are nonzero only for coordinates in which the gradient $\nabla V(\bar{y})$ is 0, i.e. $v_i \frac{\partial V}{\partial x_i}(\bar{y}) = 0$.*

Condition $(b)$ implies in particular that the Hessian $\nabla^2 V(\bar{y})$ must be rank deficient at every accumulation point, and the gradient $\nabla V(\bar{y})$ must have a zero coordinate. We stress that at this stage no assumption is made on $V$ other than its sufficient smoothness. Our result also allows the possibility of $y(t)$ having no accumulation point, in which case $(b)$ is trivially satisfied. This is for example the case for $y(t) = -t$ with $V(x) = e^x$ on $\mathbb{R}$.

The proof of Theorem 1 is presented in Section 5. We now deduce some useful special cases by strengthening some assumptions.

**Corollary 2** *Under the conditions of Theorem 1, if $\nabla^2 V$ is never rank deficient, or if its kernel never contains a vector with zero component, then exactly one of the following conditions holds:*

*(1) $y$ converges;*
*(2) $y$ admits no accumulation point. Therefore, it leaves definitively any compact set after a finite time and is unbounded.*

*In particular every bounded trajectory converges.*

**PROOF.** Under any of the additional assumptions of this corollary, no point satisfies condition $(b)$ of Theorem 1, which implies that if $y$ does not converge, it must admit no accumulation point. The rest of condition (2) follows because if there were an infinite and unbounded sequence of times at which $y$ is in a compact set $K$, it would admit an accumulation point. □

The absence of vector with zero component in the kernel of the Hessian turns out to be useful in some multi-agent systems where the Hessian can be related to a graph Laplacian, as will be seen in Section 4.2. It is also directly applicable to any function of the form $V(x) = \tilde{V}(\Pi_v x)$ where $\Pi_v$ is the orthogonal projection onto a space orthogonal to a vector $v$ with $v_i \neq 0$ for every $i$, and $\tilde{V}$ is strongly convex, i.e. to functions which are strongly convex up to an invariance in the direction $v$.

The next corollary concerns convex functions.

**Corollary 3** *Under the conditions of Theorem 1, if $\nabla^2 V(x) \succ 0 \,\forall x$ and $V$ admits a minimum $x^*$, then $y$ converges. Hence, if $V$ is strongly convex, then $y$ converges.*

**Remark 1** *The assumption of the existence of a minimum $x^*$ is needed, as $\nabla^2 V(x) \succ 0$ does not necessarily imply the existence of a minimum. See, e.g., the example $V(x) = e^x$ on $\mathbb{R}$ mentioned above.*

**PROOF.** Since $\nabla^2 V(x) \succ 0$ for every $x$, Corollary 2 implies that every bounded trajectory converges. Hence, we just need to prove that all trajectories are bounded. We suppose without loss of generality that $x^* = 0$ and $V(x^*) = 0$ and first show $V$ is radially unbounded. Indeed, if it were not, we could find arbitrary large $x$ for which $V(x) \leq M$ for some constant $M$, and convexity implies that for these $x$,

$$V\left(\frac{x}{||x||}\right) \leq \frac{||x|| - 1}{||x||} V(0) + \frac{1}{||x||} V(x) \leq \frac{M}{||x||}.$$

Hence we could find points $\frac{x}{||x||}$ of norm 1 where $V$ takes arbitrarily small value, which by the compactness of the unit sphere and the continuity of $V$ implies that $V(z) = 0$ for some $z$ of norm 1. It follows then from the convexity of $V$ that $f(\lambda) := V(\lambda z) = 0$ for all $\lambda \in [0, 1]$, and therefore that $f''(\lambda) = 0$ on $(0, 1)$. But since $f'' = z^T \nabla^2 V z$, this contradicts $\nabla^2 V(x) \succ 0$. Hence $V$ is radially unbounded, and since $V(y(t))$ is nonincreasing as $\dot{y}^T \nabla V(y(t)) \leq 0$ follows from (2), this implies the boundedness of $y$, and thus its convergence.

## 3 Tightness

We now show that the assumptions of Theorem 1 allow for situations where only condition $(b)$ holds demonstrating thus that these situations cannot be excluded without additional assumptions.

**Example 1** *Let $C = [-1, 1]^2 \subset \mathbb{R}^2$, and $V(x) = d^4(x, C)$, i.e. the fourth power of the Euclidean distance to $C$. Consider the trajectory $y(t) = (2 + e^{-t}, \sin(t))$.*

$V$ is convex, and one can verify that its Hessian is locally Lipschitz. Besides, $y$ remains in $[1, \infty) \times [-1, 1]$, on which $d(x, C) = (x_1 - 1)$ and hence $V(x) = (x_1 - 1)^4$. Furthermore the differentials in this region are

$$\nabla V(x) = \begin{pmatrix} 4(x_1 - 1)^3 \\ 0 \end{pmatrix} \quad \nabla^2 V(x) = \begin{pmatrix} 12(x_1 - 1)^2 & 0 \\ 0 & 0 \end{pmatrix}.$$

In particular $\nabla V(y(t)) = (4(e^{-t} + 1)^3, 0)^\top$, which together with $\dot{y}(t) = (-e^{-t}, \cos(t))$ implies that our assumption (2) is satisfied, and Theorem 1 applies. The trajectory does not converge, so condition $(a)$ does not hold. Its set of accumulation points $\bar{y}$ is $\{(2, c) : c \in [-1, 1]\}$, and at such points $\nabla V(\bar{y}) = (4, 0)^\top$, while the kernel of the Hessian is span$\{(0, 1)^\top\}$. These points satisfy thus condition $(b)$.

## 4 Example of application

We now demonstrate the application of our results on an example generalizing the linear platoon problem in [9], which, to the best of our knowledge, is not amenable to analysis by any available result. We show in particular how one can exploit our result in combination with additional information on the dynamics of trajectories to characterize their asymptotic behavior.

*4.1 System Description*

We consider a system of $n$ agents having each a position $y_i(t) \in \mathbb{R}$, together with a connected (undirected) graph $G$. We denote by $i \sim j$ the fact that $i$ and $j$ are neighbors on $G$, and use again $y(t)$ and $x$ respectively for the trajectory of the system and for simple vectors in $\mathbb{R}^n$.

Agents follow first-order dynamics that may be subject to unmodeled nonlinearities. They may be heterogeneous and have no common or permanent notion of time, and their progression may be slowed down or interrupted for exogenous reasons, or for example for collision avoidance. Therefore we model their dynamics as

$$\dot{y}_i(t) = h_i(t, u_i(t)), \qquad (3)$$

where $u_i(t)$ is their control input and $h_i$ an arbitrary function that preserves the sign of $u_i$ but is not known to the agents. Agents measure the relative positions of their neighbors on $G$ in an inexact manner; they have access to

$$\hat{\Delta}_{ij}(t) = y_j(t) - y_i(t) + p_{ij}(t) \quad \text{for every } j \sim i, \qquad (4)$$

for some arbitrary bounded perturbation $p_{ij}(t)$ with $|p_{ij(t)}| \leq \bar{p}_{ij}$, resulting for example from quantization, miscalibration or noise. The bounds $\bar{p}_{ij}$ are known.

We suppose that agents have certain preferences regarding their relative positions to their neighbors, characterized by "energy" functions $f_{ij}(x_j - x_i)$ satisfying the symmetry condition $f_{ij}(z) = f_{ji}(-z)$, and which they want to minimize. No assumption is made on the consistency between these functions. Agents would naturally want to minimize their total discomfort

$$V_i(x) = \sum_{j \sim i} f_{ij}(x_j - x_i),$$

with $x = (x_1, \ldots, x_n)^T$, and therefore attempt to follow $-\frac{\partial V_i}{\partial x_i} = \sum_{j \sim i} f'_{ij}(x_j - x_i)$. However, they do not have access to this exact derivative due to the measurement errors (4), and since these errors can be arbitrary, one cannot rely on them canceling out on average. Hence they will follow a robust approach and *only move* when certain that the direction is advantageous. To capture this, observe that if agent $i$ measures $\hat{\Delta}_{ij}$, (4) implies that the real $y_j - y_i$ could be anywhere in $[\hat{\Delta}_{ij} - \bar{p}_{ij}, \hat{\Delta}_{ij} + \bar{p}_{ij}]$. Hence, based on its measurements $\hat{\Delta}_{ij}$, agent $i$ can compute the following upper and lower bounds on $\frac{\partial V_i}{\partial x_i}$,

$$g_i^+(t) = \sum_{j \sim i} \sup_{p \in [-\bar{p}_{ij}, \bar{p}_{ij}]} -f'_{ij}(\hat{\Delta}_{ij}(t) + p),$$
$$g_i^-(t) = \sum_{j \sim i} \inf_{p \in [-\bar{p}_{ij}, \bar{p}_{ij}]} -f'_{ij}(\hat{\Delta}_{ij}(t) + p). \quad (5)$$

Each agent will then move only if these bounds have the same sign, as it is then sure of the sign of $\frac{\partial V_i}{\partial x_i}$. This can be achieved by

$$u_i(t) = -\max(g_i^-(t), 0) - \min(g_i^+(t), 0), \quad (6)$$

with only at most one of the terms being non-zero. The input $u_i(t)$ will thus be nonzero only if it is certain than $\frac{\partial V_i}{\partial x_i}$ is nonzero, and they will then have opposite sign. We note that the optimization problems in (5) are trivially solved by $p = \pm \bar{p}_{ij}$ when the $f_{ij}$ are convex.

*4.2 Analysis using Theorem 1*

We define the global energy function

$$V(x) := \sum_i V_i(x) = \sum_i \sum_{j \sim i} f_{ij}(x_j - x_i). \quad (7)$$

Using the symmetry of the $f_{ij}$, one can verify that, $\frac{\partial V}{\partial x_i} = 2\frac{\partial V_i}{\partial x_i}$. Moreover, our construction (3), (5) and (6) guarantee that $\dot{y}_i(t)$ could only be nonzero if one were sure it would have the opposite sign as $\frac{\partial V_i}{\partial x_i}$. Hence the trajectory of the system defined in Section 4.1 satisfies a stronger version of our main assumption (2):

$$\dot{y}_i \neq 0 \Rightarrow \dot{y}_i \frac{\partial V}{\partial x_i}|_{y(t)} < 0, \quad \forall t \geq 0, i = 1, \ldots, n. \quad (8)$$

Provided the $f_{ij}$ are twice differentiable with locally Lipschitz second derivative, we can therefore apply Theorem 1 and conclude that either $\dot{y}$ converges, or all accumulation points satisfy condition (b) of the theorem. The precise implications will depend on the nature of the functions $f_{ij}$ and of the functions $g$. We will now show when the functions $f_{ij}$ are convex and based on the non-linearity $h_i$ in (3), our result can be used to guarantee convergence to a point in the vicinity of $\arg\min V(x)$.

**Assumption 1** *Every $f_{ij}$ is locally strongly convex and admits a minimum.*
*Without loss of generality we further assume $f_{ij} \geq 0$.*

**Lemma 4** *Under Assumption 1, every trajectory of the system defined in in Section 4.1 is bounded.*

**PROOF.** We first show the existence of a (global) minimizer $x^*$ of $V$. Since every $f_{ij}$ is locally strongly convex and admits a minimum, the argument used in Corollary 3 shows that every $f_{ij}(z)$ grows unbounded when $|z|$ grows, and so does thus $\min_{i \sim j} f_{ij}(z)$. We focus now on the restriction $\tilde{V}$ of $V$ to the subspace where $\sum_i x_i = 0$. On this subspace, for any $x$ there exists a pair $k, \ell$ such that $|x_k - x_\ell| \geq ||x||_\infty$. The connectivity of $G$ implies the existence of a path of length at most $n-1$ between nodes $k$ and $\ell$, including an edge $(i, j)$ for which $|x_i - x_j| \geq \frac{1}{n-1}||x||_\infty$ Hence, it follows from the definition (7) of $V$ that for large enough $x$,

$$\tilde{V}(x) \geq f_{ij}(x_j - x_i) \geq \min_{i' \sim j'} f_{i'j'}\left(\frac{||x||_\infty}{n-1}\right),$$

using that $f_{ij}(y)$ are non-decreasing for large enough $y$. Consequently $\tilde{V}(x)$ grows to infinity when $||x||_\infty$ grows. $\tilde{V}$ is thus radially unbounded, which implies it admits a minimum $x^*$ (see e.g. [3, Proposition A.8]). Furthermore, the invariance of $V$ under addition of a a multiple of $(1, 1, \ldots, 1)^T$ implies that $x^*$ is also a minimizer of $V$.

We now fix one specific minimizer $x^*$ and consider a trajectory $y$. For every time $t$ we define the set $I(t) = \arg\max_{k=1,\ldots,n} y_k(t) - x_k^*$. This definition implies that if $i \in I(t)$, then $y_j(t) - y_i(t) \leq x_j^* - x_i^*$ for every $j$. Remember that each $f_{ij}$ is convex, and thus that each $f'_{ij}(z)$ is non-decreasing with $z$. Therefore we have

$$\frac{\partial V}{\partial x_i}|_{y(t)} = -2\sum_{j \sim i} f'_{ij}(y_j(t) - y_i(t))$$
$$\geq -2\sum_{j \sim i} f'_{ij}(x_j^* - x_i^*) = \frac{\partial V}{\partial x_i}|_{x^*} = 0,$$

since $x^*$ is a minimizer of $V$. Hence, condition (8) implies that if $\dot{y}_i \neq 0$, then $\dot{y}_i < 0$, so any agent who is reaching the maximum of $y_k(t) - x_k^*$ at a given time can only have a non-positive speed. By arguments exactly parallel to those used in [9, Theorem 4], this implies that $\max_k y_k(t) - x_k^*$ is nonincreasing. A symmetric reasoning shows then also that $\min_k y_k(t) - x_k^*$ is nondecreasing [1]. Hence every $y_i(t)$ remains bounded. □

---

[1] We note that Theorem 4 in [9] should have required condition (8) for guaranteeing the monotonicity of $\min x_i(t)$ and $\max x_i(t)$. But this condition was indeed satisfied when Theorem 4 was used in Theorem 5 of the same paper.

**Lemma 5** *let $V$ be defined as in (7). Under Assumption 1, for every $x$ there holds $\nabla^2 V(x) \succeq 0$, with $\ker \nabla^2 V = \mathrm{span}\{(1,1,\dots,1)^T\}$.*

**PROOF.** It follows from the definition of $V$ and the symmetry of the $f_{ij}$ that $(\nabla^2 V)_{ij} = -2 f''_{ij}(x_j - x_i)$ if $i \sim j$ and zero else, and that the diagonal elements are $(\nabla^2 V)_{ii} = 2\sum_{j\sim i} f''_{ij}(x_j - x_i)$. Consequently $\nabla^2 V$ can be interpreted as the Laplacian of a weighted version of the graph $G$, with the weight of each edge $(i,j)$ being $2f''_{ij}(x_j - x_i) = 2f''_{ji}(x_i - x_j)$. Moreover, these weights are all positive due to the local strong convexity of the $f_{ij}$. The result follows then from classical arguments in algebraic graph theory. □

We can now establish convergence of every trajectory.

**Corollary 6** *Under Assumption 1, every trajectory $y$ of the system defined in Section 4.1 converges.*

**PROOF.** We know that Theorem 1 applies to the trajectories $y$ of the system. And it follows from Lemma 5 that the kernel of $\nabla^2 V(x)$ is always $\mathrm{span}\{(1,1,\dots,1)^T\}$, and contains thus no nontrivial vector with a zero component. Corollary 2 implies then that every bounded trajectory converges, and thus that every trajectory converges since Proposition 4 guarantees the boundedness of all trajectories. □

There remains to characterize the potential limits of the trajectories, for which we will make an additional assumption on the function $h_i$ defining the agent dynamics in (3). The current ones allow indeed for example for $h_i \equiv 0$, in which case agents never move, and nothing can therefore be said about their limits.

**Assumption 2** *For every agent $i$, $|h_i(u_i(t), t)| \geq \alpha |u_i(t)|$ holds on an infinite number of disjoint time intervals of lengths at least $\tau$, for some $\alpha, \tau > 0$.*

It follows from (4) that the relative position $\hat\Delta_{ij}(t)$ of $j$ as measured by $i$ is at least $y_j(t) - y_i(t) - \bar p_{ij}$. Hence the convexity of $f_{ij}$ implies that $\sup_{p\in[-\bar p_{ij}, \bar p_{ij}]} -f'_{ij}(\hat\Delta_{ij}(t) - p)$ appearing in (5) is at most $-f'_{ij}(y_j(t) - y_i(t) - 2\bar p_{ij})$. Therefore, when $y(t) = x$ we have

$$g_i^+(t) \leq \tilde g_i^+(x) := -\sum_{j\sim i} f'_{ij}(x_j - x_i - 2\bar p_{ij}),$$

$$g_i^-(t) \geq \tilde g_i^-(x) := -\sum_{j\sim i} f'_{ij}(x_j - x_i + 2\bar p_{ij}).$$

Furthermore, it follows from the control law (6) that

$$u_i(t) \in [-\tilde g_i^+(y(t)), -\tilde g_i^-(y(t))], \quad \forall t, i. \qquad (9)$$

We are now ready to characterize the limiting points.

**Proposition 7** *Under Assumptions 1 and 2, every trajectory $y$ of the system defined in Section 4.1 converges. Moreover, $\lim_{t\to\infty} y(t) \in S^p := \cap_{i=1,\dots,n} S_i^p$, with $S_i^p := \{x : \tilde g_i^-(x) \leq 0 \leq \tilde g_i^+(x)\}$.*

The shape of $S^p$ depends on the properties of the $f_{ij}$, the structure of the graph $G$ and on the magnitude of the perturbation bounds $\bar p_{ij}$. Moreover, observe that it depends continuously on these bounds, and that when they are all zero, $S^p = \{x : \nabla V(x) = 0\} = \arg\min_x V(x)$. Proposition 7 proves thus convergence to a point in the vicinity of the minimum set of the potential $V$, at a distance that depends on the magnitude of the perturbations.

**PROOF.** Convergence to a point $y^*$ is guaranteed by Corollary 6. Suppose, to obtain a contradiction, that $y^* \notin S_i^p$ for some $i$, and assume without loss of generality that $\tilde g_i^+(y^*) < 0$. Then, since $\tilde g_i^+$ is continuous with respect to $x$, there exists a time $T$ after which we always have $\tilde g_i^+(y(t)) \leq \frac{1}{2}\tilde g_i^+(y^*) < 0$, and therefore $u_i(t) \geq -\frac{1}{2}\tilde g_i^+(y^*) > 0$ by (9). The sign-preserving character of $h_i$ and Assumption 2 imply then that $\dot x_i = h_i(t, u_i) \geq 0$ for all $t \geq T$, and $\dot x_i = h_i(t, u_i) \geq \alpha \frac{1}{2}\tilde g_i^+(y^*) > 0$ on an infinite number of disjoint intervals of length at least $\tau > 0$, which contradicts the convergence of $y_i(t)$ and thus of $y$. Hence we must have $y^* \in S_i^p$ for every $i$. □

## 5 Proof of Theorem 1

### 5.1 Introduction and proof structure

For the ease of reading, we will slightly abuse notations and use $\nabla V_i$ to denote $\frac{\partial V}{\partial x_i}$ and $\nabla V_i(z)$ to denote $\frac{\partial V}{\partial x_i}|_z$.

We first observe that, although we did not assume $V$ to be bounded from below, this assumption is automatically satisfied along the trajectory if there exists an accumulation point.

**Lemma 8** *Under the assumptions of Theorem 1, if $y$ admits an accumulation point $\bar y$, then $V(y(t)) \geq V(\bar y)$ for all $t$, and $\lim_{t\to\infty} V(y(t)) = V(\bar y)$.*

**PROOF.** It follows from assumption (2) that

$$\frac{d}{dt} V(y(t)) = \sum_{i=1}^n \nabla V_i(y(t)) \dot y_i(t) \leq 0,$$

implying that $V(y(t))$ is non-increasing. Since $y(t)$ gets arbitrary close to $\bar y$ for arbitrarily large times, the continuity of $V$ implies $\limsup_{t\to\infty} V(y(t)) \geq V(\bar y)$, and thanks to the monotonicity of $V(\bar y(t))$ we have then

$$\inf_t V(y(t)) = \lim_{t\to\infty} V(y(t)) = \limsup_{t\to\infty} V(y(t)) \geq V(\bar y). \qquad □$$

We now define

$$K_i := \{x : \nabla V_i(x) = 0\},$$

the set on which the $i^{th}$ coordinate of the gradient of $V$ cancels. These sets are closed by continuity of $\nabla V$. We say that an accumulation point $\bar{y}$ is *locally K-minimal* if there is a non-trivial ball centered on $\bar{y}$ containing no accumulation point that belongs to a smaller number of $K_i$ than $\bar{y}$. We first prove in Section 5.2 the result for locally K-minimal accumulation points, and will then extend it to the general case in Section 5.3 using topological arguments.

The intuition behind our proof is the following. In the non-trivial case of the theorem where the trajectory $y$ does not converge but admits a (K-minimal) accumulation point $\bar{y}$, this trajectory must repeatedly approach $\bar{y}$ and then leave it at a non-vanishing distance. We will exploit this to define a "direction" $v$ that is (asymptotically) followed infinitely often when the trajectory leaves $\tilde{y}$. We will argue that for those $i$ for which $\bar{y} \in K_i$, there must hold $(Hv)_i = 0$ with $H$ the Hessian of $V$ at $\bar{y}$, because otherwise there would be an accumulation point of the form $\bar{y} + \delta v$ at which $\nabla V_i \neq 0$, i.e. that does not belong to $K_i$, contradicting the local K-minimality of $\bar{y}$. We will also argue that for those $i$ for which $\bar{y} \notin K_i$ i.e. $\nabla V_i \neq 0$, there must hold $v_i = 0$, for otherwise, following $v$ would result in an impossible repeated decrease of energy $\nabla V(\bar{y})v \leq \nabla V_i(\bar{y})v_i < 0$ (where we use our assumption (2)). Hence there we will have $v^T H v = \sum_i v_i (Hv)_i = 0$, i.e. the direction $v$ is in the kernel of $H$. The analysis of the structure of this $v$ will then give condition (b).

5.2 *K-minimal accumulation points*

**Proposition 9** *Under the assumptions of Theorem 1, if $y$ does not converge, every locally K-minimal accumulation point $\bar{y}$ satisfies the following condition: The kernel of $\nabla^2 V(\bar{y})$ contains a nontrivial vector $v$ such that $v_i \nabla V_i(\bar{y}) = 0$ for every $i$, i.e. the entries of $v_i$ are possibly nonzero only for coordinates in which the gradient is 0.*

We suppose $y$ does not converge and fix a locally K-minimal accumulation point $\bar{y}$ (in the absence of such point, the claim trivially holds). We may re-index the coordinates without loss of generality in such a way that $\bar{y}$ belongs to $K_1, \ldots, K_k$ and not to the $n - k$ other $K_i$, with $k$ potentially equal to 0. This choice and the local K-minimality of $\bar{y}$ imply that the two following conditions are satisfied for any sufficiently small $\epsilon$, and hence we assume them to be satisfied in the sequel for the values of $\epsilon$ considered.
(i) $\bar{B}(\bar{y}, 3\epsilon) \cap K_i = \emptyset$ for $i > k$, where $\bar{B}$ denotes the closed ball.
(ii) there are no accumulation points on less than $k$ sets $K_i$ within $\bar{B}(\bar{y}, 3\epsilon)$.

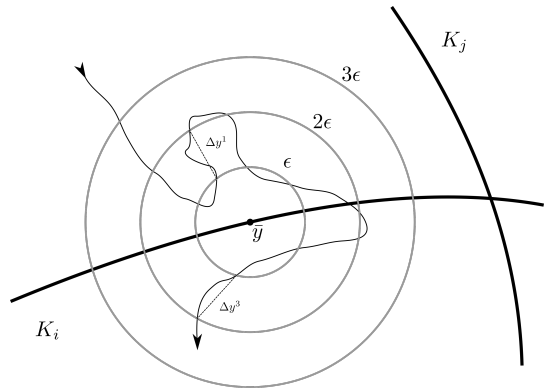

Figure 1. Representation of different constructions in the proof of Proposition 9 for a given accumulation point $\bar{y}$: (i) the sets $K_i$, with $\bar{y}$ belonging to one of them, $K_i$, so that $k = 1$; (ii) the three balls of radius $\epsilon$, $2\epsilon$ and $3\epsilon$ centered on $\bar{y}$, which do not intersect $K_j$; and (iii) the trajectory $y$ and some vectors $\Delta y^m$ connecting each time a point of the trajectory at distance $\epsilon$ from $\bar{y}$ to a subsequent point at distance $2\epsilon$ from $\bar{y}$ (these points are respectively $y(t_1^m)$ and $y(t_2^m)$, and are not represented here)

We first show that locally, the trajectory will be asymptotically constrained towards the $k$ kernel spaces $K_1, \ldots, K_k$ prescribed by $\bar{y}$. In the representation of Figure 1, this means that $y$ will approach closely $K_i$ in the long run for $i = 1, \ldots, k$.

Since $y(t)$ does not converge to $\bar{y}$, then for $\epsilon$ small enough, we can find arbitrary large $t$ for which $y(t) \notin B(\bar{y}, 2\epsilon)$. On the other hand, since $\bar{y}$ is an accumulation point, we can also find arbitrary large $t$ for which $y(t) \in B(\bar{y}, \epsilon)$. Hence, there exists a sequence of disjoint intervals $[t_1^m, t_2^m]$ with $t_1^m, t_2^m \to \infty$ such that $y(t_1^m) \in \partial B(\bar{y}, \epsilon)$, $y(t_2^m) \in \partial B(\bar{y}, 2\epsilon)$ and $y([t_1^m, t_2^m]) \in \bar{B}(\bar{y}, 2\epsilon)$. We let then $\Delta y^m = y(t_2^m) - y(t_1^m)$, be the vector linking the beginning and end of these pieces of trajectories, and already observe that $||\Delta y^m|| \in [2\epsilon, 4\epsilon]$. See again Figure 1 for a representation.

Since these vectors live in a compact set, they have an accumulation point $\Delta y^\epsilon$, still with $||\Delta y^\epsilon|| \in [2\epsilon, 4\epsilon]$. This vector can be interpreted as a direction in which the trajectory repeatedly leaves the small neighborhood of $\bar{y}$. We will show that $(\Delta y^\epsilon)^T H \Delta y^\epsilon = O(\epsilon^3)$ for $H$ the Hessian of $V$ at $\bar{y}$. For this purpose we show that $|(H \Delta y^\epsilon)_i| = O(\epsilon^2)$ for $i \leq k$ and $\Delta y_i^\epsilon = 0$ for $i > k$.

**Claim 1:** For $i \leq k$, for any increasing diverging sequence of times $t^m$ such that $y(t^m) \in \bar{B}(\bar{y}, 3\epsilon)$, there holds

$$\lim_{m \to \infty} d(y(t^m), K_i) = 0 \text{ and } \lim_{m \to \infty} \nabla V_i(y(t^m)) = 0$$

In particular $\nabla V_i(y(t_2^m)) \to 0$ and $\nabla V_i(y(t_1^m)) \to 0$.

**PROOF.** If $\lim_{m \to \infty} d(y(t^m), K_i) = 0$ did not hold, there would be an infinite subsequence $\bar{y}(t^{m'})$ at a dis-tance larger than $\delta > 0$ from $K_i$, which would admit an accumulation point $y^\delta \in \bar{B}(\bar{y}, 3\epsilon)$ since $y(t^m)$ remains in that compact set. Moreover, $y^\delta$ would not belong to any $K_j$ for $j > k$ since $\bar{B}(\bar{y}, 3\epsilon)$ does not intersect with any such $K_j$. Hence we would have an accumulation point of $y$ belonging to less than $k$ sets $K_i$ in contradiction with $\bar{y}$ being $K$-minimal. The second part of the claim follows then by continuity of the gradient [2]. □

We now show how Claim 1 implies the direction $\Delta y^\epsilon$ is "not too far" from the kernel of the first $k$ rows of the Hessian of $V$.

**Claim 2:** Let $H = \nabla^2 V(\bar{y})$. For any $i \leq k$ there holds $|(H\Delta y^\epsilon)_i| \leq C\epsilon^2$, for some $C$ possibly depending on $\bar{y}$ but not on $\epsilon$.

**PROOF.** We prove

$$\limsup_{m \to \infty} |(H\Delta y^m)_i| \leq C\epsilon^2, \quad (10)$$

which implies the result by definition of $\Delta y^\epsilon$ as an accumulation point of $\Delta y^m$. For this purpose we first show that the difference of the gradient $\nabla V(y(t_2^m)) - \nabla V(\bar{y})$ can be approximated by $H\Delta y^m$ up to $O(\epsilon^2)$. Indeed, we can write $\nabla V(y(t_2^m)) - \nabla V(\bar{y})$ as the following integral

$$\int_{s=0}^1 \nabla^2 V\left(\bar{y} + (y(t_2^m) - \bar{y})s\right)(y(t_2^m) - \bar{y})ds$$
$$= \int_{s=0}^1 H(y(t_2^m) - \bar{y})ds$$
$$+ \int_{s=0}^1 \left(\nabla^2 V\left(\bar{y} + (y(t_2^m) - \bar{y})s\right) - H\right)(y(t_2^m) - \bar{y})ds$$

Since the Hessian is assumed to be locally Lipschitz continuous, we have, for a Lipschitz constant $L(\bar{y})$,

$$\left\|\nabla^2 V\left(\bar{y} + (y(t_2^m) - \bar{y})s\right) - H\right\| \leq L(\bar{y})s \left\|y(t_2^m) - \bar{y}\right\|$$
$$= O(\epsilon),$$

with the implicit constant only depending on $\bar{y}$. Hence, slightly abusing the $O(\epsilon)$ notation for the sake of conciseness, there holds

$$\nabla V(y(t_2^m)) - \nabla V(\bar{y})$$
$$= H(y(t_2^m) - \bar{y}) + \int_{s=0}^1 O(\epsilon)(y(t_2^m) - \bar{y})ds$$
$$= H(y(t_2^m) - \bar{y}) + O(\epsilon^2),$$

---

[2] We need to work on the closed ball. Because it is compact and $\nabla V$ is continuous, Heine-Cantor theorem implies uniform continuity, hence approaching the set $K_i$ on which $\nabla V_i$ is 0 implies that $\nabla V_i$ goes to 0.

where we have used $\|y(t_2^m) - \bar{y}\| = 2\epsilon$. Similarly $\nabla V(y(t_1^m)) - \nabla V(\bar{y}) = H(y(t_1^m) - \bar{y}) + O(\epsilon^2)$. Hence

$$\nabla V(y(t_2^m)) - \nabla V(y(t_1^m)) = H\Delta y^m + O(\epsilon^2). \quad (11)$$

By Claim 1, we know that $\nabla V_i(y(t_2^m)) \to 0$ and $\nabla V_i(y(t_1^m)) \to 0$. Therefore, it follows from (11), applied to each component $i = 1, \ldots, k$, that

$$\nabla V_i(y(t_2^m)) - \nabla V_i(y(t_1^m)) = H_{i:}\Delta y^m + O(\epsilon^2) \to 0,$$

which implies (10) and thus the claim. □

Next, we show that $\Delta y_i^\epsilon = 0$ for $i > k$. The idea of the proof is that every $\Delta y_i^m$, of which $\Delta y_i^\epsilon$ is an accumulation point, results in a proportional decrease of energy "along the $i$ coordinate" that cannot be compensated by the other coordinates due to our principal condition (2).

**Claim 3:** $\Delta y_i^\epsilon = 0$ for $i > k$.

**PROOF.** We show that $\lim_{m \to \infty} \Delta y_i^m = 0$ for $i > k$, which implies the claim as $\Delta y_i^\epsilon$ is an accumulation point of $\Delta y_i^m$.

Since $\nabla V_i(\bar{y}) \neq 0$ for $i > k$ by definition of $k$, the continuity of $\nabla V$ implies that for sufficiently small $\epsilon$, we have $|\nabla V_i(x)| > c > 0$ for some $c > 0$ for all $x \in \bar{B}(\bar{y}, 2\epsilon)$. Therefore, since $y(t) \in \bar{B}(\bar{y}, 2\epsilon)$ for $t \in [t_1^m, t_2^m]$, we have

$$|\Delta y_i^m| \leq \int_{t_1^m}^{t_2^m} |\dot{y}_i|dt \leq \frac{1}{c}\int_{t_1^m}^{t_2^m} |\nabla V_i(y)| \, |\dot{y}_i|dt.$$

Our main assumption on coordinate-wise decrease (2) implies that $|\nabla V_i(y)| \, |\dot{y}_i| = -\nabla V_i(y)\dot{y}_i$, and generally that $-\nabla V_j(y)\dot{y}_j \geq 0$ for every $j$. Hence,

$$|\Delta y_i^m| \leq -\frac{1}{c}\int_{t_1^m}^{t_2^m} \left(\nabla V_i(y)\dot{y}_i + \sum_{i \neq j} \nabla V_j(y)\dot{y}_j\right) dt$$
$$= \frac{1}{c}\left(V(y(t_1^m)) - V(y(t_2^m))\right).$$

This last inequality holds for every $m$, so that

$$\sum_m |\Delta y_i^m| \leq \frac{1}{c}\sum_m (V(y(t_1^m)) - V(y(t_2^m))) < \infty,$$

as $V(y(t))$ is non-increasing and the overall decrease of $V(y(t))$ is finite by Lemma 8. Therefore, there holds $|\Delta y_i^m| \to 0$ as $m \to \infty$, which implies $\Delta y_i^\epsilon = 0$. □

**Claim 4:** Let $H = \nabla^2 V(\bar{y})$. If $\epsilon$ is small enough, there holds

$$\left|(\Delta y^\epsilon)^T H(\Delta y^\epsilon)\right| \leq C'\epsilon^3 \quad (12)$$

for some constant $C'$ depending only on $\bar{y}$.

**PROOF.** For $i \leq k$, it follows from Claim 2 that

$$|\Delta y_i^\epsilon (H \Delta y^\epsilon)_i| \leq ||\Delta y^\epsilon|| \, |(H \Delta y^\epsilon)_i| \leq 4\epsilon C \epsilon^2 =: C' \epsilon^3.$$

For $i > k$, we have $\Delta y_i^\epsilon (H \Delta y^\epsilon)_i = 0$ because $\Delta y_i^\epsilon = 0$ by Claim 3, so (12) holds. □

We are now ready to prove Proposition 9.

**PROOF.** First, remember that $||\Delta y^\epsilon|| \in [2\epsilon, 4\epsilon]$, hence any sequence of $\frac{\Delta y^\epsilon}{\epsilon}$ admits an accumulation point. In particular, among the small enough $\epsilon$ there exists a sequence of $\epsilon_\ell$ converging to 0 and a vector $v$ with $||v|| \in [2, 4]$ such that $v = \lim_{\ell \to \infty} \frac{1}{\epsilon_\ell} \Delta y^{\epsilon_\ell}$.

For $i > k$, Claim 3 implies that $\frac{1}{\epsilon_\ell} \Delta y_i^{\epsilon_\ell} = 0$ and thus $v_i = 0$. Besides, from Claim 4, we have

$$\begin{aligned}
|v^T H v| &= \lim_{\ell \to \infty} \left| \left(\frac{\Delta y^{\epsilon_\ell}}{\epsilon_\ell}\right)^T H \left(\frac{\Delta y^{\epsilon_\ell}}{\epsilon_\ell}\right) \right| \\
&= \lim_{\ell \to \infty} \frac{1}{\epsilon_\ell^2} \left| (\Delta y^{\epsilon_\ell})^T H (\Delta y^{\epsilon_\ell}) \right| \\
&\leq \lim_{\ell \to \infty} \frac{1}{\epsilon_\ell^2} C' \epsilon_\ell^3 = \lim_{\ell \to \infty} C' \epsilon_\ell = 0.
\end{aligned}$$

So we have found a nonzero $v$ in the kernel of $H$ such that $v_i = 0$ for all $i > k$, i.e. for all $i$ for which $\nabla V_i(y) \neq 0$, which establishes Proposition 9. □

5.3 *Generalization to all accumulation points*

We now prove Theorem 1 by extending the result of Proposition 9 to all accumulation points, whether K-minimal or not.

**Lemma 10** *Under the assumptions of Theorem 1, the set of $x \in \mathbb{R}^n$ satisfying property (b) in Proposition 9 is closed.*

**PROOF.** Let $T$ be the set of points $x$ satisfying the property $(b)$ in Proposition 9. Observe that $x \in T$ if and only if the following system admits a non-trivial solution

$$\begin{pmatrix} \nabla^2 V(x) \\ \text{diag}(\nabla V(x)) \end{pmatrix} v = 0,$$

i.e. if the matrix of this system is of rank smaller than $n$. Since the rank of a matrix is the size of its largest non-singular square submatrix, this condition can be checked by checking that all $n \times n$ submatrices have a zero determinant, i.e. that a finite number of continuous functions of $x$ are zero (since $\nabla^2 V(x)$ and $\nabla V(x)$ are continuous).

The set $T$ is thus an intersection of zero sets of continuous functions and is therefore closed. □

The next lemma will allow us to deduce that the set of $K$-minimal accumulation points is dense within the set of accumulation points, i.e. its closure contains all accumulation points. We state it for functions defined on subsets of $\mathbb{R}^n$, but it actually directly extends to general topological spaces, with the same proof.

**Lemma 11** *Let $S \subset \mathbb{R}^n$ and $g : S \to \mathbb{N}$. The set of local minima of $g$ is dense in $S$.*

**PROOF.** Let $M_k \subset S$ be the set of locally minimal points with value $k$. By definition,

$$M_k = g^{-1}(k) \setminus \overline{\bigcup_{i=0}^{k-1} g^{-1}(i)}. \quad (13)$$

We show by induction that $\overline{\bigcup_{i=0}^k M_i} = \overline{\bigcup_{i=0}^k g^{-1}(i)}$. For $k = 0$ this is immediate because (13) becomes $M_0 = g^{-1}(0)$. Let us now assume the relation holds for $k-1$. Using (13), we may write

$$\overline{M_k \cup \bigcup_{i=0}^{k-1} M_i} = \overline{\left(g^{-1}(k) \setminus \overline{\bigcup_{i=0}^{k-1} g^{-1}(i)}\right) \cup \overline{\bigcup_{i=0}^{k-1} g^{-1}(i)}}$$

$$= \overline{\left(g^{-1}(k) \cup \overline{\bigcup_{i=0}^{k-1} g^{-1}(i)}\right)} = \overline{\bigcup_{i=0}^k g^{-1}(i)},$$

which confirms the induction step. Since $S$ is the domain of $g$, we have then

$$S = g^{-1}(\mathbb{N}) = \bigcup_{i=0}^\infty g^{-1}(i) \subseteq \bigcup_{k=0}^\infty \overline{\bigcup_{i=0}^k g^{-1}(i)}$$

$$= \bigcup_{k=0}^\infty \overline{\bigcup_{i=0}^k M_i} \subseteq \overline{\bigcup_{i=0}^\infty M_i},$$

i.e. the closure of local minima covers $S$. □

To complete the proof of the main theorem, we let $S$ be the set of accumulation points of $y$, and define on this set the function $g$ assigning to each point the number of sets $K_i$ to which it belongs. Observe that the set $S_{\min}$ of K-minimal accumulation points is exactly the set of local minima of $g$. Hence it follows from Lemma 11 that $S_{\min}$ is dense in $S$, and thus that $S \subseteq \bar{S}_{\min}$. Now Proposition 9 states that, in the absence of convergence of $y$, every point of $S_{\min}$ satisfies condition $(b)$, and we have seen in Lemma 10 that the set of points satisfying this condition is closed. Hence every point of $S \subseteq \bar{S}_{\min}$ also satisfies that condition.

# 6 Conclusions and Open Research Directions

We have extended the results of [9] to general energy functions as opposed to quadratic positive semi-definite ones, thus significantly increasing their applicability, and clarified the possible impact of vectors with zero components in the Hessian kernel. Our results allow establishing the convergence of trajectories under very simple and easily verifiable assumptions, and guarantee in other situations simple and strong properties for the accumulation points of the trajectory.

We hope our results will serve as a useful tool for the analysis of the evolution of various systems, to take a shortcut in confirming convergence when otherwise there is a high complexity in the description of the dynamics. One may think of multi-agent interactions with communication issues, cooperation or race conditions, measurement errors and quantizations, exogenous randomness, and more.

Note that currently our results per se do not provide information on the convergence speed, but this is a consequence of an approach applicable to trajectories with potentially arbitrarily slow convergence. There remains, however, several open questions.

*Guaranteeing convergence:* We have seen in Example 1 that the alternative $(b)$ to the convergence of the trajectory cannot simply be discarded from the possible conclusions of Theorem 1. However, it might be possible to strengthen Theorem 1 by modifying our coordinate-wise decrease assumption (2). Observe indeed that Example 1 involves the trajectory freely moving along coordinates for which the corresponding gradient is zero. Hence it would not satisfy the stronger assumption (8), requiring that the derivative of the $i^{th}$ coordinate of $y$ should have opposite sign as the corresponding coordinate of $\nabla V$, and must be zero if the latter is zero. This assumption was satisfied in our example in Section 4, and exploited to prove boundedness and convergence of that system. Whether or not (8) is in general a sufficient condition for convergence is an open question. An even stronger assumption would also force $\dot{y}_i(t)$ to be nonzero when $\frac{\partial V}{\partial x_i}|_{y(t)} \neq 0$, but this would significantly decrease the applicability of the result, as it would forbid in most situations coordinates from fully stopping.

*Coordinate-free formulation:* Finally, the application of our results strongly depends on the choice of coordinates. They can be extended by embedding a change of coordinates, but formulating a truly coordinate-independent version of Theorem 1 remains an open perspective.

# Acknowledgements

We are grateful to the anonymous reviewers for their suggestions and comments. J. Hendrickx was supported by the "RevealFlight" Concerted Research Action (ARC) of the Fédération Wallonie-Bruxelles, and by the F.R.S.-FNRS via the Incentive Grant for Scientific Research (MIS) "Learning from Pairwise Comparisons" and the "KORNET" project. B. Gerencsér was supported by NRDI (National Research, Development and Innovation Office) grant KKP 137490 and the János Bolyai Research Scholarship of the Hungarian Academy of Sciences.